\documentclass[11pt,letterpaper]{amsart}

\usepackage{tikz}
\tikzstyle{every node}=[circle, draw, fill=black!50,
                        inner sep=0pt, minimum width=3pt]

\usepackage{amsfonts,amsmath,latexsym,color,epsfig,hyperref,enumitem, amssymb,bbm}

\newtheorem{theorem}{Theorem}[section]
\newtheorem{lemma}{Lemma}[section]

\newtheorem{cor}[lemma]{Corollary}
\newtheorem{prop}[lemma]{Proposition}

\renewcommand{\le}{\leqslant}
\renewcommand{\ge}{\geqslant}

\newcommand{\parag}[1]{\vspace{2mm}

\noindent{\bf #1} }

\def\qed{\ifvmode\mbox{ }\else\unskip\fi\hskip 1em plus 10fill$\Box$}

\input{epsf}

\makeatletter
\def\Ddots{\mathinner{\mkern1mu\raise\p@
\vbox{\kern7\p@\hbox{.}}\mkern2mu
\raise4\p@\hbox{.}\mkern2mu\raise7\p@\hbox{.}\mkern1mu}}
\makeatother

\def\R{\mathbb R}
\def\Z{\mathbb Z}

\def\mc{\mathcal}

\title{On sets of orthogonal exponentials on the disk}
\author{Dmitrii Zakharov}
\thanks{Zakharov's research was supported by the Jane Street Graduate Fellowship.
}
\address{Department of Mathematics, Massachusetts Institute of Technology, Cambridge, MA 02139, USA}
\email{zakhdm@mit.edu}

\date{}

\begin{document}

\begin{abstract}
We show that if $A$ is a set of mutually orthogonal exponentials with respect to the unit disk then $|A \cap [-R, R]^2| \lesssim_\varepsilon R^{3/5+\varepsilon}$ holds. This improves the previous bound of $R^{2/3}$ by Iosevich--Kolountzakis \cite{Iosevich2013}. The main new ingredient in the proof is a discretized version of Marstrand's slicing theorem.
\end{abstract}

\maketitle

\section{Introduction}

Let $D \subset \R^2$ denote the unit disk in the plane. Consider a set of points $A \subset \R^2$ such that the collection of harmonics $\{e^{2\pi i x \cdot a}\}_{a \in A}$ is pairwise orthogonal in $L^2(D)$, i.e.
\begin{equation}\label{eq1}
    \int_D e^{2\pi i x \cdot (a-a')} dx=0,
\end{equation}
for any $a \neq a' \in A$. How `large' can the set $A$ be?

This question was first raised by Fuglede \cite{Fuglede} in connection to the general question on which regions $\Omega \subset \R^n$ admit an orthogonal basis of harmonics in $L^2(\Omega)$. Fuglede conjectured that this is the case if and only if $\Omega$ tiles $\R^n$ by translations and this became a subject of extensive study in subsequent years. While the general conjecture turned out to be false in all dimensions at least $3$ \cite{ farkas2006fuglede, kolountzakis2006tiles, matolcsi2005fuglede, tao2004fuglede}, it is true e.g. when $\Omega$ is a convex set \cite{lev2022fuglede} and remains open in dimensions $d=1, 2$.  We refer to \cite{lev2019fuglede} and references therein for more details.

In the original paper, Fuglede observed that if we take $\Omega$ to be the unit disk $D$, a set which definitely does not tile the plane, then there are no such a basis and, in fact, any collection $A$ of mutually orthogonal harmonics on $D$ must be finite. Alternative proofs and extensions of this result to higher dimensions and other convex bodies with smooth boundary were obtained in \cite{Fuglede2001, iosevich1999fourier, iosevich2001convex, iosevich2003combinatorial, kolountzakis2004distance}.

The orthogonality condition (\ref{eq1}) can be rewritten as $\widehat 1_D(a-a') = 0$ for any $a\neq a' \in A$. The Fourier transform $\widehat 1_D$ is a radial function and the absolute values of its zeroes are precisely the zeros of the Bessel function $J_1(2\pi r)$. So (\ref{eq1}) becomes equivalent to
\begin{equation}\label{eq2}
    |a - a'| = r_n \text{ for some }n\ge 1\text{  and any }a \neq a'\in A, 
\end{equation}
where $r_1 < r_2 <  \ldots$ denote the zeros of $J_1(2\pi r)$. Based on the expectation that there are no non-trivial algebraic relations between zeros of the Bessel function, Fuglede \cite{Fuglede} speculated that perhaps any set $A \subset\R^2$ satisfying (\ref{eq2}) must have size at most 3 (since there always is a determinant-like algebraic relation between the six lengths among any 4 points in the plane). 
The algebraic independence of $\{r_n\}$ is strongly believed to be true but has not been established yet. 
On the other hand, in support of this prediction, Iosevich and Jaming \cite{Iosevich2008} showed that $|A| \le 3$ holds when the sequence of numbers $r_n$ is replaced by truncations of their Taylor series or similar approximations. Being unable to show any unconditional upper bound on the size of $A$, they also showed that at least it has to be a fairly sparse set: for any $R \ge 1$ we have
$$
|A \cap [-R, R]^2| \lesssim R.
$$
This bound was later refined by Iosevich and Kolountzakis \cite{Iosevich2013} to $|A \cap [-R, R]^2| \lesssim R^{2/3}$. In this note we make a further improvement:
\begin{theorem}\label{thm}
    Suppose that $A \subset \R^2$ is a set of points such that (\ref{eq2}) holds. Then for any $R \ge 1$ and $\varepsilon > 0$ we have $|A \cap [-R, R]^2| \lesssim_{\varepsilon} R^{3/5+\varepsilon}$.
\end{theorem}

Similarly to previous arguments, the proof is based on the asymptotic formula \cite{Abram}
\begin{equation}\label{eq3}
    r_n = \frac{n}{2} + \frac{1}{8} + O\left(\frac{1}{n}\right).
\end{equation}
The crucial feature of the sequence $\{r_n\}$ is that it is robustly sum-free: $|r_{n}+r_{m} - r_{k}| \ge c$ holds for any $n, m, k \ge 1$ and some absolute constant $c >0$. Our argument applies any set $A \subset \R^2$ whose distance set $\Delta(A) = \{|a-a'|, ~a\neq a' \in A\}$ is robustly sum-free and $c$-separated. It could be interesting to study how large can sets with this property (and perhaps other similar additive constraints on $\Delta(A)$) be in general. For example, is there such a set $A \subset [-R, R]^2$ of size at least $R^\epsilon$ for some constant $\epsilon>0$? 

In light of (\ref{eq3}) it is very natural to wonder what happens if we replace the sequence $\{r_n\}$ by the sequence of integers $\{n\}$. In this case one can construct fairly dense sets of integer distance sets by either putting points on a line or on a circle. In a remarkable development, Greenfeld, Iliopoulou and Peluse \cite{greenfeld2024integer} managed to show that for any integer distance set $A \subset [-R, R]^2$ there is a line or a circle $C$ such that $|A\setminus C| \le \log^{O(1)} R$. Note that it is necessary to consider the line or circle $C$ on this result: an arithmetic progression of length $R$ forms an integer distance set on a line inside $[-R,R]^2$ and there is an integer distance set of size $R^{c/\log \log R}$ on an appropriately chosen circle $C \subset [-R, R]^2$. On the other hand, if we insist that the distances in $A$ belong to the set $\{\frac{n}{2}+\frac{1}{8}\}$ (i.e. we truncate the approximation (\ref{eq3})) then $A$ has at most 2 points on any line and thus the result from \cite{greenfeld2024integer} would imply that $|A| \le R^{C/\log \log R}$.

Bounds of this form are significantly stronger bound than what we can do with our method. However the argument in \cite{greenfeld2024integer} crucially relies on the estimates on the number of low height rational points on some algebraic surfaces and so it is unclear whether it can be adapted to the case when distances lie in the set of Bessel zeros $\{r_n\}$ instead of integers $\{n\}$ or the shifted integers $\{\frac{n}{2}+\frac{1}{8}\}$.

Finally let us say a couple of words about the proof. Roughly speaking, Iosevich and Koloutzakis \cite{Iosevich2013} showed that if $A\subset \R^2$ satisfies the condition (\ref{eq2}) then it puts the following restrictions on $A$:
\begin{itemize}
    \item[(i)] If there are two points in $A$ at a small distance then $|A|$ is small,
    \item[(ii)] Any three points in $A$ with large pairwise distances cannot lie in a very thin strip,
\end{itemize}
 see Lemmas \ref{lem1} and \ref{lem2} below for precise statements. The idea now is to split into two cases depending on how $A$ `looks like' on different scales using a notion of dimension for discrete sets of points. We observe that if a set $A \subset [-R, R]^2$ is `at least $(1+\epsilon)$-dimensional' then by, the classical Marstrand slicing theorem \cite{marstrand1954some}, the intersection of $A$ with a generic thin strip has dimension at least $\epsilon$. This then leads to a contradiction with (ii) if $|A|$ is too large. Otherwise, $A$ is `at most $(1+\epsilon)$-dimensional' and this implies that $A$ contains pairs of points which are unusually close to each other, leading to some tension with (i). To make this precise we prove a discretized version of Marstrand's slicing theorem and use the notion of Katz--Tao $(\delta, s, C)$-sets, which we review in the next section. The idea of studying the `branching pattern' of a discrete set of points to find a particular subconfiguration of points also appears in the work of Cohen, Pohoata and the author \cite{cohen2023new} on the Heilbronn's triangle problem. The problem asks for the smallest area of a triangle determined by an arbitrary collection of $n$ points in the unit square. For that problem, we also study the branching pattern and use tools from the projection theory: namely, the recent radial projection theorems due to Orponen, Shmerkin and Wang \cite{orponen2024kaufman} and the high-low method.

We use standard asymptotic notation, e.g. $A\lesssim_\epsilon B$ means that $|A|\le C |B|$ for some constant $C$ depending on $\epsilon$ (and the constant is absolute if there are no subscripts), an expression $A\sim B$ means $A\lesssim B$ and $B \lesssim A$ and so on. 

\parag{Acknowledgements.} I thank Cosmin Pohoata for telling me about the problem and suggesting a connection with Heilbronn's triangle problem.
I thank Rachel Greenfeld for stimulating discussions during the Discrete Geometry workshop at Oberwolfach in January 2024 and helpful comments on an earlier version of this paper.

\section{Preliminaries}

We are going to use some basic ideas from projection theory. Throughout the proof, a $w$-square refers an axis-aligned square with side $w$. Let $Q \subset \R^2$ be a $w$-square and $A \subset \R^2$ a finite set. For some $\delta \in (0,1), s \in [0,2], C>0$, we say that $A$ is a {\em $(\delta, s, C)$-set (relative to $Q$)} if
\begin{equation}\label{eq4}
|A \cap Q'| \le C (u/w)^s |A \cap Q|,    
\end{equation}
for any subsquare $Q' \subset Q$ with side $u \in [\delta w, w]$. The notion of $(\delta, s, C)$-sets serves as a discretized analogue of sets with Hausdorff dimension $s$. Variations of this definition go back to works of Katz and Tao \cite{katz2001some} and it is a very convenient measure of how well is $A$ spread out inside of $Q$. For example, a square $\delta^{-1}\times \delta^{-1}$ grid inside of $Q$ forms a $(\delta, 2, 10)$-set, a collection of $\delta^{-1}$ uniformly spaced points on the unit circle $S^1 \subset [-1,1]^2$ forms a $(\delta, 1, C)$-set and an arbitrary collection of points in a $\delta$-neighbourhood of a point $x \in Q$ forms a $(\delta, 0, 1)$-set.  
Note that the definition is invariant upon rescaling $Q$ and so we will usually assume that $Q = [0,1]^2$ and $w=1$.

For $s, \delta > 0$ and a finite set $A \subset \R^n$ define the $(s, \delta)$-Riesz energy of $A$ by
\begin{equation}\label{eqI}
I_\delta^s(A) = \sum_{a, a' \in A} \min\{\delta^{-s}, |a-a'|^{-s}\}.   
\end{equation}

The Riesz energy is another convenient way to capture $s$-dimensional subsets in $\R^n$: note that if $A \subset [0,1]^2$ is a $(\delta, s, C)$-set then $I_\delta^{s-\epsilon}(A) \lesssim_{\epsilon, C} |A|^2$ for any $\epsilon > 0$. Indeed, expand (\ref{eqI}) and collect the pairs $a,a'\in A$ by the distance $|a-a'|$:
\[
I_{\delta}^{s-\epsilon}(A) \sim \sum_{j=0}^{\log(1/\delta)} \sum_{a, a' \in A:~ |a-a'| \sim 2^{-j}} 2^{(s-\epsilon) j}
\]
and note that for each $a \in A$ there are at most $C 2^{-sj} |A|$ choices for $a' \in A$ with $|a-a'| \sim 2^j$. Thus, 
\[
I_{\delta}^{s-\epsilon}(A) \lesssim \sum_{j=0}^{\log(1/\delta)} C 2^{-\epsilon j} |A|^2 \lesssim_{\epsilon,C} |A|^2. 
\]

Say that two $2\times \delta$ tubes $T, T'$ are essentially distinct if $|T\cap T'| \le 1.9 \delta$. Fix a collection $\mc T_0$ of $\sim \delta^{-2}$  essentially distinct $2\times \delta$ tubes such that for any line $\ell$ there is $T \in \mc T_0$ such that $\ell \cap [0,1]^2 \subset T$. More precisely, for a sequence of angles $\theta = c n \delta$, where $n=0,1, \ldots$ and $c>0$ is a small fixed constant, consider a collection $\mc T_0(\theta)$ of $\sim \delta^{-1}$ tubes inclined by the angle $\theta$ and covering the unit square. Then define $\mc T_0 = \bigcup_\theta \mc T_0(\theta)$. 

Marstrand \cite{marstrand1954some} showed that if $K \subset \R^2$ is a set of Hausdorff dimension $1+\epsilon$ then for almost all directions $\theta$, and almost all lines $\ell$ with $K\cap \ell \neq \emptyset$, the intersection $K \cap \ell$ has Hausdorff dimension $\epsilon$ (for some appropriate notions of `almost all'). The next proposition is a discretized version of this statement. Namely, we replace the set $K$ of Hausdorff dimension $1+\epsilon$ with a $(\delta, 1+\epsilon, C)$-set $A \subset [0,1]^2$, replace lines $\ell$ by $2\times \delta$ tubes $T$ and estimate the Riesz energy of the intersection $A\cap T$ instead of Hausdorff dimension of $K\cap \ell$.

\begin{prop}\label{prop}
Fix $\delta, \epsilon >0$ and $C \ge 1$. Let $A \subset [0,1]^2$ be a $(\delta, 1+\epsilon, C)$-set and $\epsilon' < \epsilon$. Let $\mc T \subset \mc T_0$ be the set of $2\times \delta$ tubes $T$ such that $|A \cap T| \in [\delta^{1+\epsilon'} |A|, \delta^{1-2\epsilon'}|A|]$ and $I_\delta^{\epsilon-\epsilon'}(A\cap T) \lesssim \delta^{-4\epsilon'} |A\cap T|^2$. Then $|\mc T| \gtrsim_{\epsilon', \epsilon, C} \delta^{-2+2\epsilon'}$.    
\end{prop}

\begin{proof}
    First, we establish that there are many $2\times \delta$ tubes $T$ with $|A\cap T| \approx \delta |A|$. Let $\pi_\theta: \R^2 \rightarrow \R$ denote the orthogonal projection in direction $\theta$. Then for any $a, a' \in A$ and $t<1$ we have
    \begin{equation}\label{eq:sum}
    \sum_{\theta} \min\{\delta^{-t}, |\pi_\theta(a) - \pi_\theta(a')|^{-t}\} \sim_t \delta^{-1} \min\{\delta^{-t}, |a-a'|^{-t}\},    
    \end{equation}
    where the sum is taken over $\theta = n c\delta$, $n=0,1, \ldots, [\frac{\pi}{c\delta}] $ and $c >0$ is a small fixed constant. So taking $t=1-\epsilon'$ and summing (\ref{eq:sum}) over all pairs $a, a' \in A$ gives
    $$
    \sum_{\theta} I^{1-\epsilon'}_\delta(\pi_\theta(A)) \sim_{\epsilon'} \delta^{-1} I^{1-\epsilon'}_\delta(A) \lesssim \delta^{-1} |A|^2,
    $$
    where the sum is over $\theta = n c\delta$, $n=0, 1, \ldots$. Thus, there is a set of $\sim \delta^{-1}$ directions $\theta$ such that $I^{1-\epsilon'}_\delta(\pi_\theta(A)) \lesssim |A|^2$. Fix one such direction $\theta$.

    On the other hand, note that for each tube $T \in \mc T_0(\theta)$ every pair $a, a' \in A\cap T$ contributes at least $\delta^{\epsilon'-1}$ to the $(1-\epsilon', \delta)$-Riesz energy of $\pi_\theta(A)$ and so we can estimate
    \begin{equation}\label{eq5}
    I_\delta^{1-\epsilon'}(\pi_\theta(A)) \gtrsim \sum_{T \in \mc T_0(\theta)} \delta^{\epsilon'-1} |A\cap T|^2.    
    \end{equation}
    By dyadic pigeonhole, we can find a dyadic integer $M \in 2^\Z$ such that the set $\mc T(\theta)$ of tubes $T \in \mc T_0(\theta)$ with $|T\cap A| \in [M \delta |A|, 2M\delta |A|)$ covers at least $\frac{c}{\log (1/\delta)}$ fraction of $A$. 
    That is, 
    $$
    \frac{|A|}{\log (1/\delta)} \lesssim \sum_{T \in \mc T(\theta)} |A\cap T| \lesssim M \delta |A| |\mc T(\theta)|
    $$
    which gives $|\mc T(\theta)| \gtrsim \delta^{-1} M^{-1} / \log(1/\delta)$ and $M \ge \delta^{\epsilon'}$ (for large enough $\delta$). Lower bounding the sum in (\ref{eq5}) by the contribution of tubes $T \in \mc T(\theta)$ then gives
    $$
    |\mc T(\theta)|  \delta^{\epsilon'-1} M^2 \delta^2 |A|^2 \lesssim  I_\delta^{1-\epsilon'}(\pi_\theta(A)) \lesssim |A|^2
    $$
    and using the lower bound on $|\mc T(\theta)|$ we obtain that $M \le \delta^{-2\epsilon'}$ holds (for large enough $\delta$). Let $\mc T_1$ be the union of sets $\mc T(\theta)$ over all directions $\theta$ for which $\pi_\theta(A)$ has $(1-\epsilon', \delta)$-Riesz energy bounded by $C|A|^2$ (by the above, there are $\sim \delta^{-1}$ such $\theta$). Thus, we constructed a collection of tubes $\mc T_1$ of size $\gtrsim \delta^{2\epsilon'-2}$ such that $|A\cap T| \in [\delta^{1+\epsilon'}|A|, \delta^{1-2\epsilon'}|A|]$ holds for any $T \in \mc T_1$.

    For a pair $a\neq a' \in A$ let $\mc T(a, a')$ be the set of tubes from $\mc T_0$ containing $a$ and $a'$. We have $|\mc T(a,a')| \sim \min\{\delta^{-1}, |a-a'|^{-1}\}$ for any $a, a' \in A$. Using this we can write
    $$
    \sum_{T \in \mc T(a, a')} \min\{\delta^{-t}, |a-a'|^{-t}\} \sim \min\{\delta^{-t-1}, |a-a'|^{-t-1}\},
    $$
    for any $t$. So summing this over all pairs $a, a' \in A$ with $t = \epsilon-\epsilon'$:
    $$
    \sum_{T \in \mc T_0} I_\delta^{\epsilon-\epsilon'}(A \cap T) \sim I_\delta^{1+\epsilon-\epsilon'}(A) \lesssim |A|^2.
    $$
    Restricting this sum to $\mc T_1$ we conclude that at least (say) half of the tubes $T \in \mc T_1$ satisfy 
    $$
    I_\delta^{\epsilon-\epsilon'}(A\cap T) \lesssim |A|^2 |\mc T_1|^{-1} \lesssim \delta^{-2\epsilon'} \delta^2 |A|^2 \lesssim \delta^{-4\epsilon'} |A\cap T|^2.
    $$
    This concludes the proof.
\end{proof}

\begin{cor}\label{lem3}
    Let $Q$ be a $w$-square and let $A \subset Q$ be a $(\delta, 1+\epsilon, C)$-set relative to $Q$, where $\epsilon, C > 0$ and $\delta$ is sufficiently small in terms of $\epsilon$ and $C$. Then there are $a_1, a_2, a_3 \in A$ such that $|a_i - a_j| \gtrsim_{\epsilon, C} \delta^{\epsilon} w$ for any $i\neq j$ and which lie in a strip of width $\delta w$.
\end{cor}

\begin{proof}
    By rescaling we may assume that $Q = [0,1]^2$ and $w=1$. By Proposition \ref{prop}, there exists a $2\times \delta$ tube $T$ such that $|A \cap T| \in [\delta^{1+\epsilon'} |A|, \delta^{1-2\epsilon'}|A|]$ and $I^{\epsilon/2}(A\cap T) \lesssim \delta^{-4\epsilon'} |A\cap T|^2$. For a ball $B$ of radius $r \ge \delta$ we have a trivial estimate
    $$
    |B \cap A\cap T|^2 r^{-\epsilon/2} \lesssim I^{\epsilon/2}_\delta(T \cap A) \lesssim \delta^{-4\epsilon'} |A\cap T|^2
    $$
    so if we take $r = \delta^{10\epsilon'/\epsilon}$ we get $|B\cap A\cap T| \lesssim \delta^{\epsilon'/2} |A \cap T|$. So we can easily find at least 3 points in $A\cap T$ with pairwise distances $\gtrsim \delta^{10\epsilon'/\epsilon}$. Take $\epsilon' = \epsilon^2/10$ and we are done.
\end{proof}

\section{Proof of Theorem \ref{thm}}

We will use the following results, obtained by Iosevich--Kolountzakis in \cite{Iosevich2013}.

\begin{lemma}[Corollary 1 \cite{Iosevich2013}]\label{lem1}
    Let $a_1, a_2,a_3 \in \R^2$ be such that $|a_i - a_j| \in  \{r_n\}$ and $|a_i - a_j| \ge L$ for any $i \neq j$. Then points $a_1, a_2, a_3$ cannot all belong to a strip of width $C L^{1/2}$, for some constant $C > 0$.
\end{lemma}

The proof of Lemma \ref{lem1} only uses the `sum-free' condition $|r_n+r_m-r_k| \ge c$ of the set $\{r_n\}$.

\begin{lemma}[Theorem 1 \cite{Iosevich2013}]\label{lem2}
    Let $A \subset \R^2$ be such that (\ref{eq2}) holds and let $t$ be the minimal distance between points in $A$. Then $|A| \le C' t$ for some constant $C'$.
\end{lemma}

The proof of Lemma \ref{lem2} is based on the Erd{\H o}s hyperbola method \cite{erdos1945integral, solymosi2003note} which was used to study integer distance sets, and also relies on Lemma \ref{lem1} and the asymptotic formula (\ref{eq3}). 

Now let $R\ge 1$ and $A \subset [-R, R]^2$ be a set such that for any $a\neq a' \in A$ we have $|a-a'| \in \{r_n\}_{n \ge 1}$. Fix some $\epsilon >0$ and $u > 0$. Let $Q \subset [-R, R]^2$ be an axis-aligned square which maximizes the expression
$$
|A\cap Q| w^{-1-\epsilon}, 
$$
where $w$ is the side of $Q$ and we have $w \in [u R, 2 R]$. We claim that $A\cap Q$ is a $(\delta, 1+\epsilon, 1)$-set relative to $Q$, where $\delta = uR/w$. Indeed, by the maximality of $Q$, for any $w'$-square $Q' \subset Q$ with $w' \in [u R, w]$ we have
$$
|A \cap Q'| w'^{-1-\epsilon} \le |A \cap Q| w^{-1-\epsilon}
$$
$$
|(A\cap Q) \cap Q'| \le (w'/w)^{1+\epsilon} |A\cap Q|,
$$
giving the claim.

First consider the case when $w \le K u R$, for some large constant $K = K(\epsilon)$. 
By applying the maximality property of $Q$ to $Q' = [-R,R]^2$ we get:
$$
|A\cap Q| w^{-1-\epsilon} \ge |A \cap [-R, R]^2| (2R)^{-1-\epsilon}, 
$$
$$
|A \cap Q| \ge (w/2R)^{1+\epsilon} |A\cap [-R, R]^2| \gtrsim u^{1+\epsilon} |A|.
$$
So if $|A| > C u^{-1-\epsilon}$ for a sufficiently large constant $C = C(\epsilon)$ then the set $A\cap Q$ contains at least two distinct points $a, a'$. We then have $|a-a'| \lesssim w \lesssim_\epsilon u R$ and so by Lemma \ref{lem2} we get $|A| \lesssim u R$. So in the first case we get $|A| \lesssim \max\{u^{-1-\epsilon}, u R\}$.

Now suppose that $w > K u R$. Then $\delta = w/uR \le 1/K$ and so Corollary \ref{lem3} applies if we take $K=K(\epsilon)$ large enough. So we obtain a triple of points $a_1, a_2, a_3 \in A$ such that $|a_i-a_j| \gtrsim L = \delta^\epsilon w$ for $i\neq j$ which lie in a strip of width $\delta w= u R$. By Lemma \ref{lem1}, we conclude that 
$$
u R \gtrsim L^{1/2} = (\delta^\epsilon w)^{1/2},
$$ 
and since $\delta=u R/w$ we get $w \lesssim (u R)^{2+2\epsilon}$. 
On the other hand, we have $|A \cap Q| \gtrsim (w/R)^{1+\epsilon} |A|$. So by covering the square $Q$ with at most $100 |A \cap Q|$ squares with side $10|A\cap Q|^{-1/2}$, we conclude that $A\cap Q$ contains a pair of points $a, a'$ at distance at most 
$$
|a-a'| \lesssim w |A\cap Q|^{-1/2} \lesssim w^{\frac{1-\epsilon}{2}} R^{\frac{1+\epsilon}{2}} |A|^{-1/2}.
$$
And so Lemma \ref{lem2} gives
$$
|A| \lesssim |a-a'| \lesssim w^{\frac{1-\epsilon}{2}} R^{\frac{1+\epsilon}{2}} |A|^{-1/2},
$$
$$
|A|^3 \lesssim w^{1-\epsilon} R^{1+\epsilon} \lesssim (u R)^2 R^{1+\epsilon},
$$
$$
|A| \lesssim u^{2/3} R^{1+\epsilon}.
$$
So by combining the two cases we obtain $|A| \lesssim \max\{u^{-1-\epsilon}, u R, u^{2/3} R^{1+\epsilon}\}$. Taking $u = R^{-3/5}$ gives the desired bound $|A| \lesssim R^{3/5+\epsilon}$.

\parag{Remark.} One can check that `the worst' set $A$ for our argument has the following branching pattern: it is $1$-dimensional between scales $R^{4/5}$ and $R$, 2-dimensional between scales $R^{3/5}$ and $R^{4/5}$ and 0-dimensional below scale $R^{3/5}$ (i.e. any two points of $A$ are at least $R^{3/5}$ distance apart). In particular, if one wants to improve the exponent $3/5$ further then one may essentially assume that $A$ has this special form.

\bibliographystyle{amsplain0.bst}
\bibliography{main}

\end{document}